\documentclass[11pt]{article}
\usepackage{xypic}
\usepackage{amssymb}
\usepackage{amsfonts}
\usepackage{amsmath}
\usepackage[mathscr]{eucal}
\usepackage{latexsym}
\usepackage[frenchb]{babel}

\newtheorem{teor}{Th{\'e}or{\`e}me}
\newtheorem{prop}[teor]{Proposition}
\newtheorem{cor}[teor]{Corollaire}                                             
\newtheorem{rem}[teor]{Remarque}    
\newtheorem{lema}[teor]{Lemme}

\begin{document}
\author{Gentiana Danila}
\title{Sections de la puissance tensorielle  du fibr{\'e} tautologique  sur le sch{\'e}ma de Hilbert ponctuel d'une surface}
\date{\ \ \ \ }
\maketitle

\def\vs{\vskip 0.5cm}

\def\has{{\rm H}}
\def\te{{\rm T}}
\def\be{{\rm B}}
\def\de{{\rm D}}
\def\tZ{\tilde{Z}}
\def\tmu{\tilde{\mu}}
\def\ttmu{\tilde{\tmu}}
\def\tpi{\tilde{\pi}}
\def\tnu{\tilde{\nu}}
\def\tp{\tilde{p}}
\def\barpi{\bar{\pi}}
\def\barp{\bar{p}}

\def\ka{{\rm K}}
\def\es{{\rm S}}
\def\er{{\rm R}}
\def\supp{{\rm supp}\,}
\def\esk{\es^k}
\def\esnmk{\es^{n-k}}

\def\Hom{{\rm Hom}}
\def\Pic{{\rm Pic}\,}
\def\Tor{{\rm Tor}\,}
\def\Spec{{\rm Spec}\,}
\def\Ker{{\rm Ker}\,}
\def\Im{{\rm Im}\,}
\def\uTor{\underline{{\rm Tor}}\,}
\def\pgl{{\rm PGL}\,}
\def\gl{{\rm GL}\,}


\def\maH{{\mathcal{H}}}
\def\maD{{\mathcal{D}}}
\def\maO{{\mathcal{O}}}
\def\maL{{\mathcal{L}}}
\def\maI{{\mathcal{I}}}
\def\sigm{{\mathfrak{S}}}
\def\sigmd{{\mathfrak{S}}_2}
\def\sigmn{{\mathfrak{S}}_n}
\def\sigmk{{\mathfrak{S}}_k}
\def\sigmkp{{\mathfrak{S}}_{k+1}}
\def\emfrac{{\mathfrak{m}}}

\def\pe{{\cal P}}
\def\perond{\underline{\pe}}
\def\peront{\perond\ \widetilde{}}
\def\peunu{\pe_1}
\def\pel{\pe_l}
\def\pem{\pe_m}
\def\deltap{\Delta_{\perond}}
\def\arond{a_{\perond}\;}
\def\aronds{a^*_{\perond}\;}
\def\hasrond{h_{\perond}\;}
\def\hasronds{h_{\perond}^*\;}


\def\proj{{\mathbb{P}}}
\def\pp{\proj_2}

\def\zeteka{(Z/T)^k}
\def\hilx{X^{\mbox{}^{[n]}}}
\def\tihilx{X_{\sim}^{\mbox{}^{[n]}}}
\def\trhilx{X_{(3)}^{\mbox{}^{[n]}}}
\def\ddhilx{X_{(2,2)}^{\mbox{}^{[n]}}}
\def\dhilx{\partial X^{\mbox{}^{[n]}}}
\def\hilxs{X_*^{\mbox{}^{[n]}}}
\def\hilxss{X_{**}^{\mbox{}^{[n]}}}
\def\dhilxs{\partial X_*^{\mbox{}^{[n]}}}
\def\hild{X^{\mbox{}^{[2]}}}
\def\ppd{\proj_2^{\mbox{}^{[2]}}}
\def\dhild{\partial X^{\mbox{}^{[2]}}}
\def\hilds{X_*^{\mbox{}^{[2]}}}
\def\dhilds{\partial X_*^{\mbox{}^{[2]}}}
\def\hilxnu{X^{\mbox{}^{[n,1]}}}
\def\hilxnd{X^{\mbox{}^{[n,2]}}}
\def\hilxndp{X^{\mbox{}^{[n,2]\prime}}}
\def\hilxnk{X^{\mbox{}^{[n,k]}}}
\def\hilk{X^{\mbox{}^{[k]}}}
\def\benk{\be^{\mbox{}^{[n,k]}}}

\def\bes{B_*}
\def\bess{B_{**}}
\def\xgen{X^k_{{\it gen}}}
\def\tess{T_{**}}
\def\tes{T_{*}}
\def\espe{\es_{\perond}}
\def\espess{\es_{\perond\;**}}
\def\esgen{\es_{{\it gen}}}
\def\esgenbar{\overline{\es}_{{\it gen}}}
\def\ceppl{C_{\perond,\pel}}
\def\hgen{h_{{\it gen}}\;}
\def\hgens{h^*_{{\it gen}}\;}
\def\pgen{p_{{\it gen}}\;}
\def\pgens{p^*_{{\it gen}}\;}

\def\hilt{X^{\mbox{}^{[3]}}}
\def\hilts{X_*^{\mbox{}^{[3]}}}
\def\hiltss{X_{**}^{\mbox{}^{[3]}}}
\def\dhilt{\partial X^{\mbox{}^{[3]}}}
\def\dhilts{\partial X_*^{\mbox{}^{[3]}}}
\def\hiltc{X_{3C}^{\mbox{}^{[3]}}}
\def\hiltd{X^{\mbox{}^{[3,2]}}}
\def\xd{X^2}
\def\esxn{\es^nX}
\def\esxd{\es^2X}
\def\desxd{\partial\esxd}
\def\xn{X^n}
\def\xk{X^k}
\def\xns{X^n_*}
\def\esdoi{\es^2}
\def\bedoi{\be^2}
\def\bedois{\be^2_*}
\def\benij{\be^n_{ij}}
\def\benijk{\be^n_{ijk}}
\def\benijkl{\be^n_{ij,kl}}
\def\ben{\be^n}
\def\den{\de^n}
\def\tiben{\be^n_{\sim}}
\def\tiden{\de^n_{\sim}}

\def\trben{\be^n_{(3)}}
\def\ddben{\be^n_{(2,2)}}
\def\tisigma{\Sigma_{\sim}}
\def\bens{\be^n_*}
\def\benss{\be^n_{**}}
\def\benud{\be^n_{12}}
\def\bends{\be^{\mbox{}^{[n,2]}}_*}
\def\bendij{\be^{\mbox{}^{[n,2]}}_{ij}}
\def\xij{\Xi_{ij}}
\def\xii{\Xi_{i}}
\def\xiud{\Xi_{[12]}}
\def\xiudi{\Xi_{[12]i}}
\def\xiund{\Xi_{12}}
\def\iunk{i_1,\cdots,i_k}
\def\xiunk{\Xi_{i_1\cdots i_k}}
\def\xiuk{\Xi_{i_1\cdots i_k}}

\def\deltaij{\Delta_{ij}}
\def\eij{E_{ij}}

\def\eskaxi{\esk_{\hilx}(\Xi)}
\def\xikah{\Xi^k_{\hilx}}
\def\zed{{\mathbb Z}}
\def\comp{{\mathbb C}}


\def\dedoi{{\maD_2}}
\def\dna{{\maD_n^A}}
\def\dda{{\maD_2^A}}
\def\dta{{\maD_3^A}}
\def\dnl{{\maD_n^L}}
\def\ddl{{\maD_2^L}}
\def\dkl{{\maD_k^L}}

\def\dedoid{{\maD_2^{\mbox{}^{[2]}}}}
\def\dedoit{{\maD_2^{\mbox{}^{[3]}}}}
\def\dedoin{{\maD_2^{\mbox{}^{[n]}}}}
\def\dkn{{\maD_k^{\mbox{}^{[n]}}}}
\def\dnk{{\maD_k^{\mbox{}^{[n]}}}}

\def\elen{L^{\mbox{}^{[n]}}}
\def\elden{L^{{2\mbox{}^{[n]}}}}
\def\eled{L^{\mbox{}^{[2]}}}
\def\eltr{L^{\mbox{}^{[3]}}}
\def\eld{L^2}
\def\tield{\widetilde{L^2}}
\def\eldtr{L^{{2\mbox{}^{[3]}}}}
\def\eldoi{\maL_2}
\def\elded{L^{{2\mbox{}^{[2]}}}}
\def\elij{L_{ij}}

\def\eskal{\esk\elen}


\def\pnu{p_{n1}}
\def\pinu{\pi_{n1}}
\def\pnd{p_{n2}}
\def\pind{\pi_{n2}}
\def\pnus{p_{n1}^*}
\def\pinus{\pi_{n1*}}
\def\pnds{p_{n2}^*}
\def\pinds{\pi_{n2_*}}
\def\pink{\pi_{nk}}
\def\penk{p_{nk}}
\def\pinks{\pi_{nk*}}
\def\penks{p_{nk}^*}

\def\surto{\twoheadrightarrow}
\def\ra{\rightarrow}
\def\lra{\longrightarrow}


\def\tens{\otimes}


\def\fs{{faisceau }}
\def\fx{{faisceaux }}
\def\alg{{alg{\'e}brique }}
\def\algs{{alg{\'e}briques }}
\def\th{{th{\'e}or{\`e}me }}

\def\thm{{th{\'e}or{\`e}me }}

\def\iso{{isomorphisme }}
\def\rep{{repr{\'e}sentation }}
\def\reps{{repr{\'e}sentations }}
\def\irr{{irr{\'e}ductible }}
\def\irrs{{irr{\'e}ductibles }}
\def\fib{{fibr{\'e} }}
\def\fibs{{fibr{\'e}s }}
\def\mor{{morphisme }}
\def\mors{{morphismes }}
\def\sur{{surjectif }}
\def\co{{coh{\'e}rent }}
\def\cohs{{coh{\'e}rents }}
\def\app{{application }}
\def\apps{{applications }}

\def\resp{{respectivement }}
\def\inve{{inversible }}
\def\inves{{inversibles }}

\def\sct{{section }}
\def\scts{{sections }}
\def\cano{{canonique }}
\def\canos{{canoniques }}
\def\lin{{lin{\'e}aire }}
\def\lins{{lin{\'e}aires }}
\def\schil{{sch{\'e}ma de Hilbert }}

{\bf Abstract~:} {\small We compute the space of global sections  for the tensor power of the 
tautological bundle  on the punctual
Hilbert scheme $\hilx$ of a smooth projective surface $X$ on
$\comp$. }



Running heads: Fibr{\'e} tautologique sur le sch{\'e}ma de Hilbert d'une surface

\vs Soit $X$ une surface complexe projective et lisse et $L$
un \fib \inve sur $X$. Soit $n$ un entier positif. On note $T=\hilx$ le sch{\'e}ma de Hilbert qui param{\`e}tre les sous-sch{\'e}mas de $X$ de longueur $n$. Il est lisse et projectif de dimension $2n$ (\cite{Fogarty}).

On consid{\`e}re la vari{\'e}t{\'e} d'incidence $Z\subset T\times X$ des points $(\xi,x)$ qui v{\'e}rifient $x\in \supp \xi$.

On note $\pi, p$ les projections
\begin{equation}
\label{ecuatia00}
{\diagram
Z\rto^{p}\dto_{\pi}&X\\
T.&
\enddiagram}
\end{equation}

On d{\'e}finit $\elen=\pi_*(p^* L)$. C'est un faisceau localement libre de rang $n$ sur $T$. 

On consid{\`e}re le produit fibr{\'e} $\zeteka$ au-dessus de $T$. On note encore $\pi, p$ les projections
\begin{equation}
\label{ecuatia01}
{\diagram
\zeteka\rto^{p}\dto_{\pi}&X^k\\
T.&
\enddiagram}
\end{equation}
Le \mor  $\pi:Z\to T$ {\'e}tant fini et plat, la formule de changement de base montre que
$$\pi_*p^*(L^{\boxtimes k})=(\elen)^{\otimes k}.$$
On obtient le \mor 
\begin{equation}
\label{ecuatia02}
\has^0(X,L)^{\otimes k}=\has^0(X^k,L^{\boxtimes k})\stackrel{p^*}{\to}\has^0(\zeteka,p^*L^{\boxtimes k})=\has^0(T,\pi_*p^*(L^{\boxtimes k}))=\has^0(T,(\elen)^{\otimes k}).
\end{equation}
Le groupe sym{\'e}trique $\sigmk$ agit sur les sch{\'e}mas $X^k, \zeteka$ et sur les \fx $L^{\boxtimes k},(\elen)^{\otimes k}$. Le \mor (\ref{ecuatia02}) est {\'e}quivariant.

Le but de cet article est de d{\'e}montrer:
\begin{teor} 
\label{th1}
Si $n\ge k\ge 0$ alors le \mor (\ref{ecuatia02})
$$p^*:\has^0(X,L)^{\otimes k}\to\has^0(\hilx,(\elen)^{\otimes k})$$
est un \iso de $\sigmk$-repr{\'e}sentations.
\end{teor}

\begin{cor}
\label{cor2}
Si $n\ge k\ge 0$ alors
$$\has^0(\hilx,\eskal)\simeq\esk\has^0(X,L).$$
\end{cor}

\par {\bf Preuve du corollaire \ref{cor2}:}

On consid{\`e}re les invariants par l'action de $\sigmk$ dans l'isomorphisme du \th \ref{th1}. $\Box$

\vskip 0.3cm

La remarque \ref{rema10} discute la relation entre le corollaire \ref{cor2} et le calcul de l'espace des sections globales du \fib d{\'e}terminant de Donaldson sur l'espace de modules des \fx semi-stables sur le plan projectif (la formule de Verlinde pour le plan projectif).

\vskip 0.3cm

L'id{\'e}e de la d{\'e}monstration du \th \ref{th1} consiste dans l'{\'e}tude des intersections des strates de $\zeteka$ obtenues de la stratification canonique de $X^k$. Plus pr{\'e}cis{\'e}ment, on consid{\`e}re les partitions $\perond=\{\peunu,\cdots,\pem\}$ de l'ensemble $\{1,\cdots,k\}$ en parties disjointes:
$$\{1,\cdots,k\}=\bigcup_{l=1,\cdots,m}\pel.$$
On appelle $m$ la longueur de la partition. La vari{\'e}t{\'e} $X^k$ est la r{\'e}union disjointe des ensembles localement ferm{\'e}s
$$\deltap=\{(x_1,\cdots,x_k)\ |\ x_i=x_j {\ \rm s'il\ existe \ }l{\rm \ avec\ } i,j\in\pel\}$$
param{\'e}tr{\'e}s par les partitions $\perond$ de l'ensemble $\{1,\cdots,k\}$.
On note $gen$ la partition g{\'e}\-n{\'e}\-rique \linebreak[4] $\{\{1\},\cdots,\{k\}\}$. Il lui correspond l'ouvert $\xgen\subset X^k$ des points {\`a} coordonn{\'e}es distinctes. Pour {\'e}tendre la stratification de $X^k$ {\`a} $\zeteka$ on utilise:

\begin{lema}
\label{lema3}
Le sch{\'e}ma $\zeteka$ est r{\'e}duit.
\end{lema}

\par {\bf Preuve:}

Le \mor $\pi:Z\to T$ est fini et plat. Le sch{\'e}ma $Z$ est de Cohen-Macaulay (\cite{Eisenbud}, cor. 18.17). 
De m{\^e}me le sch{\'e}ma $\zeteka$ est fini et plat au-dessus du sch{\'e}ma lisse $ T$  donc il est de Cohen-Macaulay. Le \mor $\pi:Z\to T$ est lisse au-dessus de l'ouvert $\tess$ des sch{\'e}mas qui ont pour support $n$ points distincts. Par cons{\'e}quent le sch{\'e}ma $\zeteka$ est lisse en dehors du ferm{\'e} $\pi^{-1}(T\setminus \tess)$ de codimension $1$. Il est donc r{\'e}duit.  (\cite{Mats}, p.183)$\Box$

\vskip 0.3cm

Par le lemme \ref{lema3} on trouve que $\zeteka$ co{\"\i}ncide avec le ferm{\'e} des points $(\xi,(x_1,\cdots,x_k))\in T\times X^k$ tels que $x_i\in\supp\xi, \forall i.$ Il r{\'e}sulte que $\zeteka$ est la r{\'e}union des ensembles localement ferm{\'e}s
$$\espe=\{(\xi,(x_1,\cdots,x_k))\in T\times\deltap, x_i\in\supp\xi, \forall i\}.$$
On consid{\`e}re la strate $\esgen$, non vide d'apr{\`e}s l'hypoth{\`e}se $n\ge k$. On note $\esgenbar$ sa fermeture dans $\zeteka$. Les \mors
$$\esgenbar\stackrel{i}{\to}\zeteka\stackrel{p}{\to}X^k$$
induisent les \mors
\begin{equation}
\label{ecuatia03}
\has^0(X^k,L^{\boxtimes k}) \stackrel{p^*}{\to}\has^0(\zeteka,p^*L^{\boxtimes k})\stackrel{i^*}{\to}\has^0(\esgenbar,i^*p^*L^{\boxtimes k}).
\end{equation}

Le \th \ref{th1} r{\'e}sulte des deux propositions suivantes:

\begin{prop}
\label{propo4}
Le \mor $i^*p^*$ est un isomorphisme.
\end{prop}

\begin{prop}
\label{propo5}
Le \mor $i^*$ est injectif.
\end{prop}

La d{\'e}monstration de ces propositions n{\'e}cessite un changement de base du sch{\'e}ma de Hilbert $T$ {\`a} la vari{\'e}t{\'e} $\bes$ suivante.
Soit $\xns$ l'ouvert des points $(x_1,\cdots,x_n)\in\xn$ avec au plus deux coordonn{\'e}es $x_i,x_j$ {\'e}gales. On note $D$ la r{\'e}union  des diagonales disjointes $D_{ij}=\{(x_1,\cdots,x_n)\ |\ x_i=x_j\}$ de $\xns$. On note $\bes$ l'{\'e}clat{\'e} de $D$ dans $\xns$. Le diviseur exceptionnel $E$ est la r{\'e}union disjointe des diviseurs $E_{ij}$ correspondant {\`a} chaque $D_{ij}$. Le groupe sym{\'e}trique $\sigmn$  agit sur $\bes$. Le quotient $\bes/\sigmn$ s'identifie {\`a} l'ouvert $\tes\subset T$ des sch{\'e}mas avec au plus un point double. On note $\rho:\bes\to\xns$ le \mor d'{\'e}clatement et $q:\bes\to\tes$ le \mor quotient.

On consid{\`e}re une partition $\perond=\{\peunu,\cdots,\pem\}$ de $\{1,\cdots,k\}$. 
On choisit un ordre $(\peunu,\cdots,\pem)$, et 
pour chaque $i\in\{1,\cdots,k\}$, on note $\ell(i)$ l'unique indice tel que $i\in\pe_{\ell(i)}$. On d{\'e}finit l'application $\arond:\xn\to\xk$:
$$\arond(x_1,\cdots,x_n)=(x_{\ell(1)},\cdots,x_{\ell(k)}).$$
(exemple: lorsque $k=7, \peunu=\{3,4\}, \pe_2=\{5\},\pe_3=\{1,6,7\},\pe_4=\{2\},$
$$\arond(x_1,\cdots,x_n)=(x_3,x_4,x_1,x_1,x_2,x_3,x_3))$$
On d{\'e}finit l'application associ{\'e}e {\`a} la partition ordonn{\'e}e $\perond$:
$$\hasrond=(q,\arond\circ\rho):\bes\to T\times\xk.$$
\begin{lema}
\label{lema6}
L'image de l'application $\hasrond$ est la fermeture $\overline{\es}_{\perond}$ de la strate $\espe$. L'application $\hasrond$ induit le \mor injectif:
\begin{equation}
\label{ecuatia04}
\has^0(\overline{\es}_{\perond},p^*L^{\boxtimes k})\stackrel{\hasronds}{\lra}\has^0(\bes,\hasronds p^*L^{\boxtimes k})=\has^0(X,L^{|\peunu|})\otimes\cdots\otimes\has^0(X,L^{|\pem|}).
\end{equation}
\end{lema}

\par {\bf Preuve:}

On rappelle que $\tess\subset T$ est l'ouvert des sch{\'e}mas {\`a} support $n$ points distincts. On note $\bess=q^{-1}(\tess)$ et $\espess\subset\espe$ l'ouvert $\espe\cap\pi^{-1}(\tess)$. Par d{\'e}finition, $\hasrond(\bess)=\espess$ et cela suffit pour la premi{\`e}re affirmation.
Par d{\'e}finition on a 
$\arond\circ\rho=p\circ \hasrond:\bes\to \xk.$ Par suite 
$$\displaylines{
\has^0(\bes,\hasronds p^* L^{\boxtimes k})=\has^0(\bes,\rho^*\aronds L^{\boxtimes k}) =\has^0(\bes,\rho^*(L^{|\peunu|}\boxtimes\cdots\boxtimes L^{|\pem|}\boxtimes\maO_{X^{n-m}})).}$$
Le \mor $\rho:\bes\to\xns$ est l'{\'e}clatement d'une sous-vari{\'e}t{\'e} lisse, donc $\rho_*\maO_{\bes}=\maO_{\xns}$. D'apr{\`e}s la formule de projection on obtient:
$$\displaylines{
\has^0(\bes,\rho^*(L^{|\peunu|}\boxtimes\cdots\boxtimes L^{|\pem|}\boxtimes\maO_{X^{n-m}}))=\has^0(\xns,L^{|\peunu|}\boxtimes\cdots\boxtimes L^{|\pem|}\boxtimes\maO_{X^{n-m}})=\hfill\cr
\hfill=\has^0(\xn,L^{|\peunu|}\boxtimes\cdots\boxtimes L^{|\pem|}\boxtimes\maO_{X^{n-m}})=\has^0(X,L^{|\peunu|})\otimes\cdots\otimes\has^0(X,L^{|\pem|}).}$$
On a utilis{\'e} le fait que $\xns$ est un ouvert dont le compl{\'e}mentaire est de codimension $\ge 2$ dans $\xn$ et la formule de K{\"u}nneth. $\Box$

\vskip 0.3cm

\par {\bf Preuve de la proposition \ref{propo4}:}

On consid{\`e}re la suite de \mors
$$ \bes\stackrel{\hgen}{\lra}\esgenbar\stackrel{\pgen}{\lra}\xk.$$
Elle induit la suite de \mors
\begin{equation}
\label{ecuatia05}
\has^0(X,L)^{\otimes k}=\has^0(\xk,L^{\boxtimes k})\stackrel{\pgens}{\lra}\has^0(\esgenbar,\pgens L^{\boxtimes k})\stackrel{\hgens}{\lra}\has^0(\bes,\hgens\pgens L^{\boxtimes k})
\end{equation}
o{\`u} $\pgens$ est le \mor $i^*p^*$ de l'{\'e}nonc{\'e}.
D'apr{\`e}s le lemme \ref{lema6} on a
$$ \has^0(\bes,\hgens\pgens L^{\boxtimes k})=\has^0(X,L)^{\otimes k}.$$
Il r{\'e}sulte de la d{\'e}finition que la composition de la suite (\ref{ecuatia05}) est l'identit{\'e} sur $\has^0(X,L)^{\otimes k}$. Par cons{\'e}quent $\hgens$ est surjectif. Par le lemme \ref{lema6} le \mor $\hgens$ est injectif, donc $\hgens$ est un isomorphisme et $\pgens=i^*p^*$ aussi. $\Box$

\vskip 0.3cm

\par {\bf Preuve de la proposition \ref{propo5}:}

Soit $s$ une section globale du \fs $p^*(L^{\boxtimes k})$ sur $\zeteka$, qui s'annule sur $\esgenbar$. Puisque le sch{\'e}ma $\zeteka$ est r{\'e}duit, il suffit de d{\'e}montrer que la restriction de la section $s$ {\`a} chacune des strates $\espe$ est nulle. On va d{\'e}montrer le r{\'e}sultat par r{\'e}currence descendente selon la longueur de la partition $\perond$. Le pas initial $\perond= gen$ est tautologique.

Soit $\perond=\{\peunu,\cdots,\pem\}$ une partition 
diff{\'e}rente de $gen$. Il existe un ensemble $\pel$ avec $|\pel|>1$. On consid{\`e}re une {\'e}criture de $\pel$ comme r{\'e}union disjointe $\pel=\pe'_l\cup\pe''_l$ d'ensemble non vides. On note $\peront$ la partition $\{\peunu,\cdots,\pe_{l-1},\pe'_l,\pe''_l,\cdots,\pem\}.$ C'est une partition de longueur $m+1$.
On avait not{\'e} $\tes\subset T$ l'ouvert des sch{\'e}mas avec au plus un point double. On note $\partial\tes\subset\tes$ le ferm{\'e} des sch{\'e}mas avec exactement un point double. On introduit l'ensemble $\ceppl\subset T\times\xk$:
$$\ceppl=\{(\xi,(x_1,\cdots,x_k))|\ \xi\in\partial\tes {\ \rm de \ point \ singulier \ } y,(x_1,\cdots,x_k)\in\deltap,x_i\in\supp\xi,y=x_j, \forall j\in\pel\}.$$
On a $\ceppl\subset\espe$. La proposition \ref{propo5} r{\'e}sulte des deux lemmes suivants:

\begin{lema}
\label{lema7}
L'inclusion $j:\ceppl\to\espe$ induit le \mor injectif
$$\has^0(\espe,p^*L^{\boxtimes k})\to\has^0(\ceppl,j^*p^*L^{\boxtimes k}).$$
\end{lema}

\begin{lema}
\label{lema8}
L'ensemble $\ceppl$ est inclus dans l'adh{\'e}rence $\overline{\es}_{\peront}$ de la strate $\es_{\peront}$.
\end{lema}

Effectivement, si la restriction de la section globale $s$ {\`a} la strate $\es_{\peront}$ est nulle, par le lemme \ref{lema8} sa restriction {\`a} l'ensemble $\ceppl$ est nulle et par le lemme \ref{lema7} sa restriction {\`a} la strate $\espe$ est nulle. Ainsi on a montr{\'e} le pas de r{\'e}currence de $\peront$ {\`a} $\perond$. $\Box$

\vskip 0.3cm

\par {\bf Preuve du lemme \ref{lema7}:}

Soit $i,j$ deux indices distincts dans $\pel$, et $\eij\subset\bes$ le diviseur exceptionnel. On consid{\`e}re un ordre sur la partition $\perond$ et la restriction du \mor $\hasrond$ associ{\'e} {\`a} cet ordre 
au diviseur $\eij$. De mani{\`e}re analogue {\`a} la preuve du lemme \ref{lema6} on d{\'e}montre que $\hasrond(E_{ij})=\ceppl$ et que le \mor $\hasrond$ induit le \mor injectif:
\begin{equation}
\label{ecuatia06}
\has^0(\ceppl,p^*L^{\boxtimes k})\stackrel{\hasronds}{\lra}\has^0(E_{ij},\hasronds p^*L^{\boxtimes k})\simeq\has^0(X,L^{|\peunu|})\otimes\cdots\otimes\has^0(X,L^{|\pem|}).
\end{equation}
Les \mors (\ref{ecuatia04}) et (\ref{ecuatia06}) rentrent dans le diagramme commutatif:
\begin{equation}
\label{ecuatia07}
{\diagram
\has^0(\espe,p^*L^{\boxtimes k})\rto^{\hasronds}\dto&\has^0(\bes,\hasronds p^*L^{\boxtimes k})\dto^{\simeq}\\
\has^0(\ceppl,p^*L^{\boxtimes k})\rto^{\hasronds}&\has^0(E_{ij},\hasronds p^*L^{\boxtimes k}).
\enddiagram}
\end{equation}
Les \mors horizontaux sont injectifs, d'o{\`u} la conclusion. $\Box$

\vskip 0.3cm

\par {\bf Preuve du lemme \ref{lema8}:}

Soit $(\xi,(x_1^0,\cdots, x_k^0))$ un point de $\ceppl$. Pour $\pe_s$ partie de la partition $\perond$, on note $x_{\pe_s}$ la valeur $x_i^0$ pour tout $i\in\pe_s$ (elles sont {\'e}gales). Le support du sch{\'e}ma $\xi$ est par d{\'e}finition
$$x_{\peunu}+\cdots+x_{\pe_{l-1}}+2x_{\pel}+x_{\pe_{l+1}}+\cdots+x_{\pem}+y_1+\cdots+y_{n-m+1}.$$
Il existe un germe de courbe lisse $C_t$ sur $\tes$ tel que $C_0=\xi$:
$$C_t=x_{\peunu}+\cdots+x_{\pe_{l-1}}+x'_t+x''_t+x_{\pe_{l+1}}+\cdots+x_{\pem}+y_1+\cdots+y_{n-m+1}.$$
On construit le germe de courbe lisse $D_t=(C_t,(x_1^t,\cdots,x_k^t))$ sur $\tes\times\xk$ ainsi:
\begin{itemize}
\item $x_i^t=x_i^0$ si $i\in\pe_s, s\ne l$
\item $x_i^t=x'_t$ si $i\in\pe'_l$
\item $x_i^t=x''_t$ si $i\in\pe''_l$.
\end{itemize}
Le germe $D_t$ est dans $\es_{\peront}$ pour $t\ne 0$ et il est le point initial pour $t=0$. $\Box$

\begin{rem}
\label{rema9}
{\rm On a utilis{\'e} l'hypoth{\`e}se $n\ge k$ pour pouvoir r{\'e}duire par le proc{\'e}d{\'e} des lemmes \ref{lema7} et \ref{lema8} toutes les strates $\espe\subset\zeteka$ {\`a} la strate g{\'e}n{\'e}rique $\esgen=\es_{\{\{1\},\cdots,\{k\}\}}$. Cet argument
ne fonctionne pas quand $n<k$. Par exemple pour $n=2, k=4$ les strates ``les plus g{\'e}n{\'e}rales'' $\es_{\{\{1,2\},\{3,4\}\}}$ et $\es_{\{\{1,2,3\},\{4\}\}}$ ne sont pas r{\'e}ductibles l'une {\`a} l'autre par ce proc{\'e}d{\'e}.
Lorsque $n=1,k>1$ on a par d{\'e}finition
$$\has^0(\hilx,(\elen)^{\otimes k})=\has^0(X,L^k)$$
diff{\'e}rent de $\has^0(X,L)^{\otimes k}$.}
\end{rem}

\begin{rem}
\label{rema10}
{\rm Soit $A$ un \fib \inve sur $X$. Soit $\esxn$ le quotient de la vari{\'e}t{\'e} $\xn$ par l'action du groupe sym{\'e}trique $\sigmn$. Soit $\dna=(A^{\boxtimes n})^{\sigmn}$ le \fib des invariants du \fib $A^{\boxtimes n}$. C'est un \fib \inve sur $\esxn$. On note toujours $\dna$ l'image r{\'e}ciproque de $\dna$  par le \mor de Hilbert-Chow $HC:\hilx\to\esxn$. On 
a r{\'e}duit dans \cite{Danila1} le calcul de l'espace de sections du fibr{\'e} d{\'e}terminant de Donaldson sur l'espace de modules de faisceaux semi-stables de rang $2$ sur le plan projectif au calcul des groupes de cohomologie:
$$\has^*(\hilx,\es^k\elen\otimes\dna).$$
L'argument de cet article ne s'{\'e}tend pas au calcul de $\has^0(\hilx,\eskal\otimes\dna)$. Pour sim\-pli\-fi\-er l'e\-xemple on prend $n=k=2$ et $\perond=\{\{1,2\}\}, \pel=\{1,2\}$ dans le lemme \ref{lema7}. Le \mor vertical du diagramme (\ref{ecuatia07}) est:
$$\has^0(X,L^2\otimes A)\otimes\has^0(X,A)\to\has^0(X,L^2\otimes A^2)$$
qui n'est pas un \iso pour $A$ g{\'e}n{\'e}ral.}

\end{rem}

{\bf Note:} Cet article est une g{\'e}n{\'e}ralisation de l'approche de J. Le Potier qui a simplifi{\'e} de fa{\c c}on substantielle mon article original \cite{Danila3}.

{\bf Remerciements:} Je remercie M. Brion pour avoir attir{\'e} mon
attention sur le lien entre les composantes irr{\'e}ductibles du sch{\'e}ma
sym{\'e}trique $\esk(Z/T)=\zeteka/\sigmk$  et les calculs de
\cite{Danila1}. Je remercie J. Le Potier pour la simplification
apport{\'e}e au \cite{Danila3} et pour nos discussions
fructueuses. Durant les diff\'erentes \'etapes qui ont men\'e \`a cette
version je me suis r{\'e}jouie de l'ambiance amicale de l'Institut de Math{\'e}matiques de l'Universit{\'e} de Warwick.

\bigskip 

Gentiana Danila\\
Institut de Math{\'e}matiques de Jussieu, \\
Universit{\'e} Paris 7, Case Postale 7012\\
175, rue du Chevaleret, 75013 Paris \\
e-mail: gentiana@math.jussieu.fr

\end{document}